\documentclass[envcountsame]{llncs}
\usepackage{amsmath,amsfonts,amssymb}
\usepackage{array}
\usepackage[all]{xy}
\usepackage{colordvi}
\usepackage[dvips]{rotating} 
\usepackage[dvips]{graphics} 

\renewcommand{\blacksquare}{{\vrule height0.70ex width0.25em}}

\newcommand{\CC}{{\bbbc}}
\newcommand{\QQ}{{\bbbq}}
\newcommand{\PP}{{\mathbb P}}
\newcommand{\Me}{{\mathcal M}}

\newcommand{\GCD}{{\rm GCD}}
\newcommand{\GL}{{\rm GL}}
\newcommand{\SL}{{\rm SL}}
\newcommand{\Sym}{{\rm Sym}}
\newcommand{\ord}{{\rm ord}}
\newcommand{\Res}{{\rm Res}}

\newcommand{\ignore}[1]{}

\title{Classification of Genus 3 Curves\\
in Special Strata of the Moduli Space}
\author{Martine Girard, David R. Kohel}
\authorrunning{M.~Girard, D.~R.~Kohel}
\institute{
  School of Mathematics and Statistics, \\
  The University of Sydney \\
  \email{\{girard,kohel\}@maths.usyd.edu.au}
}

\begin{document}
\maketitle

\begin{abstract}
We describe the invariants of plane quartic curves --- nonhyperelliptic 
genus 3 curves in their canonical model --- as determined by Dixmier 
and Ohno, with application to the classification of curves with given 
structure.
In particular, we determine modular equations for the strata in the 
moduli space $\Me_3$ of plane quartics which have at least seven 
hyperflexes, and obtain an computational characterization of curves in 
these strata.
\end{abstract}

\section{Introduction}

The classification of curves of genus 0, 1, and 2 is aided by use of 
various geometric and arithmetic invariants.  In this work we consider 
nonhyperelliptic genus 3 curves, for which the canonical model is an 
embedding as a projective plane quartic.
The work of Hess~\cite{Hess2002,Hess2004} gives generic algorithms for 
determining the locus of Weierstrass points and for finding whether two 
curves are isomorphic.
Such a generic approach to isomorphism testing works well for curves over 
finite fields, where a small degree splitting field for the Weierstrass 
places exists, and when one wants to test only two curves.  
In this work, we investigate the geometric invariants of nonhyperelliptic 
genus 3 curves, which are much more suited to classifying curves which 
are already given in terms of their canonical embeddings. 

In particular, plane quartic curves admit explicit formulas for the 
Weierstrass locus, invariants of Dixmier and Ohno by which the curves 
may be classified up to isomorphism over an algebraically closed field, 
and moreover can be classified into strata following Vermeulen's 
characterization in terms of the number and configuration of Weierstrass 
points of weight two. 

In the generic case, Harris~\cite{Harris1979} proved that a generic 
curve of any genus over a field of characteristic zero, is expected to 
have generic Galois action on the Weierstrass points.  Thus in order 
to establish an isomorphism between the sets of Weierstrass points one 
needs in general an excessively large degree extension to apply the 
algorithm of Hess.  Thus it becomes essential to exploit any special 
structure of the Weierstrass points to facilitate this algorithm.
In this article we focus on curves whose moduli lie in special strata 
of the moduli space of genus three curves.  We use a classification
by invariants to reduce to a trivial calculation of invariants on 
certain strata of Vermeulen of dimensions 0 and 1.

\section{The Weierstrass Locus of Quartic Curves}

A nonhyperelliptic curve $C$ of genus 3 can be define via the canonical 
embedding by a quartic equation $F(X,Y,Z) = 0$ in the projective plane.
The Hessian $H(X,Y,Z)$ of the form $F(X,Y,Z)$ is defined by 
$$
H = 
\left|
\begin{array}{ccc}
\displaystyle
\frac{\partial^2 F}{\partial X^2} & 
\displaystyle
\frac{\partial^2 F}{\partial X \partial Y} & 
\displaystyle
\frac{\partial^2 F}{\partial X \partial Z} \\ \\
\displaystyle
\frac{\partial^2 F}{\partial X \partial Y} & 
\displaystyle
\frac{\partial^2 F}{\partial Y^2} & 
\displaystyle
\frac{\partial^2 F}{\partial Y \partial Z} \\ \\
\displaystyle
\frac{\partial^2 F}{\partial X \partial Z} & 
\displaystyle
\frac{\partial^2 F}{\partial Y \partial Z} &
\displaystyle
\frac{\partial^2 F}{\partial Z^2}
\end{array}
\right|\cdot
$$
This form is a sextic, which meets the curve $C$ in the 24 inflection 
points of the curve (counting multiplicities).   These inflection 
are also the Weierstrass points, hence they may be determined in an 
elementary way.  The inflection points which meet the Hessian with 
multiplicity 2 are those inflection points which meet their tangent 
line to multiplicity 4, and are called {\em hyperflexes}. 

The hyperflexes are intrinsic points of the genus 3 curve, since they 
are Weierstrass of weight~2 -- those which have a deficit of 2 in 
their gap sequences.  Thus they are preserved by any isomorphism of 
curves, and reflect the underlying geometry of the curves (rather 
than solely of a choice of projective embedding).

We focus in this article on the classification of those curves which 
have an exceptional number of hyperflexes.  The partioning of the 
Weierstrass points into those of weight 1 and weight 2 can also be 
efficiently determined since it is the singular subscheme of the 
intersection $F = H = 0$, defined by the vanishing of the Jacobian
minors:
$$
\frac{\partial F}{\partial X}\frac{\partial H}{\partial Y} - 
\frac{\partial F}{\partial Y}\frac{\partial H}{\partial X} = 
\frac{\partial F}{\partial X}\frac{\partial H}{\partial Z} - 
\frac{\partial F}{\partial Z}\frac{\partial H}{\partial X} = 
\frac{\partial F}{\partial Y}\frac{\partial H}{\partial Z} - 
\frac{\partial F}{\partial Z}\frac{\partial H}{\partial Y} = 0.
$$
The calculation of the hyperflex locus can be reduced to polynomial 
factorization, without the need for Gr{\"o}bner basis calculations.
Let $R$ be the resultant $\Res(H,F,Z)$ of degree 24 and let set 
$G(X,Y) = \GCD(R,R_X,R_Y)$.  Then $G$ determines the $(X,Y)$-coordinates 
of the hyperflex locus for which $XYZ \ne 0$.

Since plane quartics are canonical embeddings of a genus~3 curve, 
any isomorphism of such curves is induced by a linear isomorphism 
of their ambient projective planes.  As a result, the problem of 
determining isomorphisms is reduced to the intersection of a linear 
algebra problem of finding such an isomorphism and a combinatorial  
one, of mapping a finite set of Weierstrass points to Weierstrass 
points.  By combining classification of quartics by their moduli
invariants into strata determined by the numbers and configurations 
of hyperflexes, we facilitate the problem of establishing 
isomorphisms between curves.

\section{Quartic Invariants}

The $j$-invariant of an elliptic curve or the Igusa invariants of a genus~2 
curve provide invariants by which a curve of genus~1 or~2 can be classified 
up to isomorphism.  Recall that the $j$-invariant of an elliptic curve is 
defined in terms of weighted projective invariants $E_4$, $E_6$, and $\Delta$ 
such that 
$$
j = \frac{E_4^3}{\Delta} \mbox{ and } E_4^3-E_6^2 = 12^3 \Delta.
$$
Similarly Igusa~\cite{Igusa1960} defined weighted invarinants $J_2$, $J_4$, 
$J_6$, $J_8$, and $J_{10}$ with one relation $J_2J_6 - J_4^2 = 4J_8$, from 
which one can determine a set of absolute invariants (see also
Mestre~\cite{Mestre1991}).  

For genus $g \ge 2$, the moduli space of curves of genus~$g$ is a space 
of dimension $3g-3$, thus the determination of generators for the ring 
of projective invariants becomes increasingly difficult.  
However, for genus~3, Dixmier~\cite{Dixmier1987} provided an explicit 
set of 7 weighted invariants and proves their algebraic independence 
over $\CC$.  The determination of these invariants builds on on the 
explicit 19th century methods of Salmon~\cite{Salmon1879}).  
By comparison of the Poincar{\'e} series of this ring with that computed 
by Shioda~\cite{Shioda1967}, Dixmier finds that his invariants $I_3$, 
$I_6$, $I_9$, $I_{12}$, $I_{15}$, $I_{18}$, and $I_{27}$, determine a 
subring over which the ring ${\mathcal A}$ of all invariants is free 
of rank $50$.

Recently, Ohno~\cite{Ohno200x} determined a complete set of generators 
and relations for the full ring of invariants of plane quartics.   
In particular he shows that there exist six additional invariants $J_9$, 
$J_{12}$, $J_{15}$, $J_{18}$, $I_{21}$, $J_{21}$, which generate 
${\mathcal A}$ as a finite algebraic extension of 
$\CC[I_3,I_6,I_9,I_{12},I_{15},I_{18},I_{27}]$.

\subsection*{Definition of covariants and contravariants}

In this section we recall the definitions of covariants and contravariants,
and basic constructions appearing in Dixmier~\cite{Dixmier1987}, 
Ohno~\cite{Ohno200x}, and Salmon~\cite{Salmon1879}.  We first introduce 
the definitions of covariant and contravariant, following the modern 
terminology used in Poonen, Schaefer, and Stoll~\cite[\S7.1]{PSS2005}.
For other modern treatments of the subject, see 
Sturmfels~\cite{Sturmfels1993} and Olver~\cite{Olver1999}. 

Let $V = \CC^n$ be be equipped with the standard left action of $\GL_n(\CC)$, 
which induces a right action on the algebra $\CC[x_1,\dots,x_n] = \Sym(V^*)$. 
For $\gamma$ in $\GL_n(\CC)$ and $F \in \CC[x_1,\dots,x_n]$ we define this 
action by $F^\gamma(x) = F(\gamma(x))$ for all $x$ in $V$.  
We denote $\CC[x_1,\dots,x_n]_d$ the $d$-th graded components of polynomials 
homogeneous of degree $d$.

\begin{definition}
A {\it covariant} of degree $r$ and order $m$ is a $\CC$-linear function
$$
\psi: \CC[x_1,\dots,x_n]_d \longrightarrow \CC[x_1,\dots,x_n]_m,
$$
which satisfies
\begin{enumerate}
\item
$\SL_n(\CC)$-module homomorphism, i.e. 
$\psi(F^\gamma) = \psi(F)^\gamma$ for all $\gamma \in \SL_n(\CC)$, 
\item
the coefficients of $\psi(F)$ depend polynomially in the coefficients 
of $x_1^{i_1}\cdots x^n_{i_n}$, 
\item
and 
$\psi(\lambda F) = \lambda^r \psi(F)$ for all $\lambda \in \CC$. 
\end{enumerate}
\end{definition}
We note in particular that the last two conditions imply that $\psi$ 
is homogeneous of degree $r$ in the coefficients of the degree $d$ 
form $F$.  An {\it invariant} is a covariant of order $0$.

\noindent
{\bf N.B.} One ``usually'' defines a covariant to satisfy
$$
\psi(F^\gamma(\bar{x})) = \det(\gamma)^k \psi(F(x))
\mbox{ where } \gamma(\bar{x}) = x, 
$$
or equivalently,
$
\psi(F^\gamma) = \det(\gamma)^k \psi(F)^\gamma,
$
for all $\gamma$ in $\GL_n(V)$ and $x$ in $V$.  One defines $k$ to be 
the {\it weight} (or {\it index}) of $\psi$.  Clearly we then have 
the relation $2k = dr - m$.  Applying a scalar matrix $\gamma$ to an 
invariant implies that $k \equiv 0 \bmod n$.  The definition of 
Poonen, Schaefer, and Stoll admits the possibility of covariants with 
a multiplicative character.  Following the classical definitions we 
include the stronger condition above in our definition of covariants.

In order to define a contravariant, we set $\CC[u_1,\dots,u_n] = \Sym(V)$,
where $\{u_1,\dots,u_n\}$ is a basis for $V$ dual to the basis 
$\{x_1,\dots,x_n\}$ of $V^*$.
Then $\GL_n(\CC)$ has the right contragradient action on polynomials in
$\CC[u_1,\dots,u_n]$, which we denote $G^{\gamma_*}$, where $\gamma_*$ 
is the inverse transpose of $\gamma$.

\begin{definition}
A {\it contravariant} of degree $r$ and order $m$ is a $\CC$-linear function
$$
\psi: \CC[x_1,\dots,x_n]_d \longrightarrow \CC[u_1,\dots,u_n]_m,
$$
which satisfies
\begin{enumerate}
\item
$\SL_n(\CC)$-module homomorphism, i.e. 
$\psi(F^\gamma) = \psi(F)^{\gamma_*}$ for all $\gamma \in \SL_n(\CC)$, 
\item
the coefficients of $\psi(F)$ depend polynomially in the coefficients 
of $x_1^{i_1}\cdots x^n_{i_n}$, 
\item
and 
$\psi(\lambda F) = \lambda^{-r} \psi(F)$ for all $\lambda \in \CC$. 
\end{enumerate}
\end{definition}

\noindent
{\bf N.B.} As noted in~\cite{PSS2005}, we may formally identify 
$u_1,\dots,u_n$ with $x_1,\dots,x_n$, via the isomorphism 
$V \rightarrow V^*$ implied by the choice of basis for $V$.  
We nevertheless distinguish the $GL_n(\CC)$-modules structures by 
denoting the action by $G \mapsto G^{\gamma_*}$ for $\gamma$ 
in $\GL_n(\CC)$.  In our mathematical exposition we preserve the 
notational distinction between $x_i$ and $u_i$.

\subsection*{Covariant and contravariant operations}

We extend the linear pairing $V \times V^* \rightarrow \CC$ given by 
$(u_i,x_j) \mapsto \delta_{ij}$ to a differential operation
$$
D: \CC[u_1,\dots,u_n] \times \CC[x_1,\dots,x_n] \rightarrow \CC[x_1,\dots,x_n], 
$$
by identifying a monomial $u_1^{i_1} \cdots u_n^{i_n}$ of total degree 
$m$ with the operator
$$
\frac{\partial^m}{\partial^{i_1} x_i \cdots \partial^{i_n} x_n}\cdot
$$
We denote the $D(\psi,\varphi)$ by $D_\psi(\varphi)$.  By symmetry, we 
define a differential operator 
$$
D: \CC[x_1,\dots,x_n] \times \CC[u_1,\dots,u_n] \rightarrow \CC[u_1,\dots,u_n], 
$$
and denote $D(\varphi,\psi)$ by $D_\varphi(\psi)$. (We resolve the 
notational ambiguity from the arguments.)  We recall as a lemma a 
classical result.

\begin{lemma}
\label{DOperation}
Let $\varphi$ be a covariant and $\psi$ be a contravariant on 
$\CC[x_1,\dots,x_n]_d$.  Then $D_\varphi(\psi)$ is a contravariant 
of order $\ord(\psi) - \ord(\varphi)$ and and $D_\psi(\varphi)$ a 
covariant of order $\ord(\psi) - \ord(\varphi)$, both of degree 
$\deg(\varphi) + \deg(\psi)$. 
\end{lemma}

When specializing to ternary forms, we denote $(x_1,x_2,x_3)$ by 
$(x,y,z)$.  For ternary quadratic forms, Dixmier~\cite{Dixmier1987} 
used the additional operations.  
Let $\varphi$ be a ternary quadratic form in $x,y,z$ and 
$$
D(\varphi) = \frac{1}{2}\left[
\begin{array}{ccc}
\displaystyle\frac{\partial^2\varphi}{\partial x^2} & 
\displaystyle\frac{\partial^2\varphi}{\partial x\partial y} & 
\displaystyle\frac{\partial^2\varphi}{\partial x\partial z} \\ 
\displaystyle\frac{\partial^2\varphi}{\partial x\partial y} & 
\displaystyle\frac{\partial^2\varphi}{\partial y^2} & 
\displaystyle\frac{\partial^2\varphi}{\partial y\partial z} \\ 
\displaystyle\frac{\partial^2\varphi}{\partial x\partial z} & 
\displaystyle\frac{\partial^2\varphi}{\partial y\partial z} &
\displaystyle\frac{\partial^2\varphi}{\partial z^2}
\end{array}
\right]
$$
and let $D(\varphi)^*$ be its classical adjoint.  Then for $\varphi$
and $\psi$ covariant and contravariant forms, respectively, we define
$$
J_{11}(\varphi,\psi) = \langle D(\varphi),D(\psi)\rangle \mbox{ and }
J_{22}(\varphi,\psi) = \langle D(\varphi)^*,D(\psi)^*\rangle,
$$
where $\langle{A,B}\rangle$ is a matrix dot product, and 
$$
J_{30}(\varphi,\psi) = J_{30}(\varphi) = \det(D(\varphi)) \mbox{ and }
J_{03}(\varphi,\psi) = J_{03}(\psi) = \det(D(\psi)).
$$
The expressions $J_{ij}$ play a role in invariant theory of ternary 
quadratic forms, but more generally we have the following classical 
lemma. 
\begin{lemma}
Let $\varphi$ be a covariant and $\psi$ be a contravariant on 
$\CC[x,y,z]_d$, each of order $2$.  Then $J_{ij}(\varphi,\psi)$ 
is an invariant on $\CC[x,y,z]_d$ of degree $i\deg(\varphi) + 
j\deg(\psi)$.
\end{lemma}
In particular we will apply this to describe the construction of the 
complete invariants of ternary quartics by Dixmier~\cite{Dixmier1987}
and Ohno~\cite{Ohno200x}.

\vspace{1mm}

Finally, for two binary forms $F(x,y)$ and $G(x,y)$ of degrees $r$ and $s$, 
we define the $k$-th {\it transvectant} of is defined to be  
$$
(F,G)^k = \frac{(r-k)!(s-k)!}{r!s!}\Big(
\frac{\partial^2}{\partial x_1\partial y_2} - 
\frac{\partial^2}{\partial y_1\partial x_2}\Big)^k
F(x_1,y_1)G(x_2,y_2) \Big|_{(x_i,y_i)=(x,y)}
$$

\begin{lemma}
\label{BinaryQuartic}
Let $F(x,y) = a_{40}x^4 + 4a_{31}x^3y + 6a_{22}x^2y^2 + 4a_{13}xy^3 + a_{04}y^4$ 
be a binary quartic form, and set $G = (F,F)^2$.  Then $F$ has invariants 
$\sigma$ and $\psi$ defined by
$$
\begin{array}{l}
\displaystyle
\sigma = \frac{1}{2}(F,F)^4 = a_{40}a_{04} - 4a_{31}a_{13} + 3a_{22}^2, \mbox{ and } \\
\displaystyle
\psi = \frac{1}{6}(F,G)^4 = 
a_{40}a_{22}a_{04} - a_{40}a_{13}^2 - a_{31}^2a_{04} + 2a_{31}a_{22}a_{13} - a_{22}^3.
\end{array}
$$
The invariant $\sigma^3-27\psi^2$ is the discriminant of the form $F(x,y)$ (up to a 
scalar).
\end{lemma}

\subsection*{Covariants and contavariants of quartics}

In this section we use the above construction to define the invariants 
of Dixmier and Ohno classifying ternary quartics, i.e.~genus~3 
curves of general type.  As noted above, we denote $(x_1,x_2,x_3)$ by 
$(x,y,z)$, and dual $(u_1,u_2,u_3)$ by $(u,v,w)$.  The polynomial rings 
$k[x,y,z]$ and $k[u,v,w]$ represent coordinate rings of the ambient 
projective space $\PP^2$ and the dual projective space $(\PP^2)^*$, 
respectively, of the quartic $F(x,y,z) = 0$.  Where necessary, 
$k[x_1,x_2,x_3,u_1,u_2,u_3]$ will be the bi-graded coordinate ring 
of $\PP^2 \times (\PP^2)^*$.

\vspace{2mm}
\noindent
{\bf N.B.} In the spirit of the classical literature, we speak of 
covariants and contravariants of a quartic form $F(x,y,z)$, though 
formally the covariant or contravariant is a function from $\CC[x,y,z]_4$ 
to $\CC[x,y,z]_m$ or $\CC[u,v,w]_m$.  Similarly, we may express a 
homogeneous form as 
$$
F(x_1,\dots,x_n) = \sum_{i_1,\dots,i_n} 
  \frac{d!}{i_1!\dots i_n!} a_{(i_1,\dots,i_n)} x_1^{i_1}\cdots x_n^{i_n},
$$
where the sum is over all indices with $i_1+\cdots+i_n = d$.  The 
calculation of invariants is thus normalized to be primitive with respect 
to such {\it classically integral} forms (as in Lemma~\ref{BinaryQuartic}).  
In the case of quartics, the constructions and expressions often require 
the primes 2 and 3 to be invertible, even if the final invariants can be 
made well-defined in these characteristics.  In what follows we follow 
Dixmier and Ohno in normalizing the expressions to be primitive with 
respect to the coefficients $\{a_{ijk}\}$, and only at the end provide 
the scalars by which the invariants must be normalized to be integral 
with respect to the coefficients $a_{ijk}$ of a integral form
\begin{equation}
\label{IntegralQuartic}
F(x,y,z) = \sum_{i,j,k} a_{i,j,k} x^iy^jz^k.
\end{equation}
In the spirit of Igusa's article on genus 2 curves~\cite{Igusa1960}, the 
determination of a complete set of integral invariants over any ring, 
and their algorithmic construction, remains open.

The first covariants at our disposal for a form $F$ is the form itself 
(i.e. the identity covariant) and the Hessian $H$.  We additionally 
require two contravariants $\sigma$ and $\psi$, from which the Dixmier 
and Ohno invariants are derived by the operations of the previous 
section.  The contravariant $\sigma$ appears in Salmon~\cite{Salmon1879}
(\S92 and \S292), has degree 2 and order 4, and the construction of 
the degree 3 and order 6 contravariant $\psi$ appears (in Salmon 
\S92,~p.78).
Formally intersect $ux + vy + wz = 0$ with the form $F$, and setting
$w = 1$, eliminate $z$ a binary quartic $R(x,y) = F(x,y,-ux-vy)$. 
Then the invariants $\sigma$ and $\psi$ of Lemma~\ref{BinaryQuartic},
rehomogenized with respect to $w$, provide us with the covariants 
$\sigma(u,v,w)$ and $\psi(u,v,w)$.

We can now define a system of covariants and contravariants for ternary 
quartics, from the covariants $F$ and $H$ and contravariants $\sigma$ 
and $\psi$.
\begin{center}
\begin{tabular}{lll}
\setlength{\itemsep}{20mm}
\setlength{\baselineskip}{20mm}
Covariants & \quad\quad & Contravariants \\ \hline
$\displaystyle\tau = {12}^{-1} D_{\rho}(F)$ & &
$\displaystyle \rho = {144}^{-1} D_{F}(\psi)$   \\
$\displaystyle
\xi = {72}^{-1} D_{\sigma}(H)$ & &
$\displaystyle \eta = {12}^{-1} D_{\xi}(\sigma)$ \\
$\displaystyle
\nu = {8}^{-1} D_{\eta}D_{\rho}(H)$ & &
$\displaystyle \chi = {8}^{-1} D_{\tau}^2(\psi)$
\end{tabular}
\end{center}
Subsequently we can define the invariants of Dixmier:
$$
\begin{array}{lll}
\displaystyle I_3 = {144}^{-1} D_\sigma(F),\quad\quad & 
\displaystyle I_9 = J_{11}(\tau,\rho),\quad\quad\quad & 
\displaystyle I_{15} = J_{30}(\tau),\\
\displaystyle I_6 = {4608}^{-1}(D_\psi(H) - 8I_3^2),\quad\quad &
\displaystyle I_{12} = J_{03}(\rho),\quad &
\displaystyle I_{18} = J_{22}(\tau,\rho).
\end{array}
$$
together with the discriminant $I_{27}$.  Dixmier~\cite{Dixmier1987} 
proved that these invariants are algebraically independent over $\CC$ 
and generate a subring of the ring ${\mathcal A}$ of ternary quartic 
invariants of index 50.  Ohno~\cite{Ohno200x} proved computationally
that the additional six invariants 
$$
\begin{array}{lll}
\displaystyle J_9 = J_{11}(\xi,\rho),\quad &
\displaystyle J_{15} = J_{30}(\xi),\quad &
\displaystyle I_{21} = J_{03}(\eta), \\
\displaystyle J_{12} = J_{11}(\tau,\eta),\quad\quad   &
\displaystyle J_{18} = J_{22}(\xi,\rho),\quad\quad  &
\displaystyle J_{21} = J_{11}(\nu,\eta),
\end{array}
$$
generate ${\mathcal A}$; he moreover determined a complete set of 
algebraic relations for the ring ${ \mathcal A} = \CC[I_k,J_l]$.

In the following table we summarise the covariant and contravariant 
degrees and orders, as can be determined from Lemma~\ref{DOperation},
beginning with the forms $F$, $H$, $\sigma$ and $\psi$.
\begin{center}
\begin{tabular}{c|ccccccccc}
\multicolumn{8}{c}{Covariants} \\
       & \multicolumn{6}{c}{$\ord$} \\
$\deg$ &\quad 0\quad & \quad 1\quad  & \quad 2 \quad & 
\quad 3 \quad & \quad 4 \quad & \quad 5 \quad & \quad 6 \\ \hline
   1   &   &   &   &   &\quad F  &   &   \\
   2   &   &   &   &   &   &   &   \\
   3  &$I_3$&   &   &   &   &   & \quad H \\
   4   &   &   &   &   &   &   &   \\
   5   &   &   &\quad $\tau$, $\xi$ &&&& \\
   6  &$I_6$&   &   &   &   &   &   \\
   7   &   &   &   &   &   &   &   \\
   8   &   &   &   &   &   &   &   \\
   9 &$I_9,J_9$&   &   &   &   &   &   \\
  10   &   &   &   &   &   &   &   \\
  11   &   &   &   &   &   &   &   \\
  12&$I_{12},J_{12}$& &   &   &   &   &   \\
  13   &   &   &   &   &   &   &   \\
  14   &   &\quad $\nu$&   &   &   &   &   \\
\end{tabular}
\quad
\begin{tabular}{c|ccccccccc}
\multicolumn{8}{c}{Contravariants} \\
       & \multicolumn{6}{c}{$\ord$} \\
$\deg$ & \quad 1\quad  & \quad 2 \quad & 
\quad 3 \quad & \quad 4 \quad & \quad 5 \quad & \quad 6 \\ \hline
   1   &   &   &   &   &   &   \\
   2   &   &   &   &\quad $\sigma$&&  \\
   3   &   &   &   &   &   &\quad $\psi$ \\
   4   &   &\quad $\rho$&&   &   &   \\
   5   &   &   &   &   &   &   \\
   6   &   &   &   &   &   &   \\
   7   &  &\quad $\eta$& &   &   &   \\
   8   &   &   &   &   &   &   \\
   9   &   &   &   &   &   &   \\
  10   &   &   &   &   &   &   &   \\
  11   &   &   &   &   &   &   &   \\
  12   &   &   &   &   &   &   &   \\
  13   &&\quad $\chi$&&  &   &   &   &   \\
  14   &   &   &   &   &   &   &   \\
\end{tabular}
\end{center}
As noted above, the natural normalization for the invariant to be integral 
depends whether one considers classically integral forms or integral forms.  
On integral forms one normalizes the Dixmier--Ohno invariants as follows:
$$
\begin{array}{c}
2^4 3^2I_3,\quad
2^{12} 3^6 I_6,\quad
2^{12} 3^8 I_9,\quad
2^{16} 3^{12} I_{12},\quad
2^{23} 3^{15} I_{15},\quad
2^{27} 3^{17} I_{18},\quad
2^{40} I_{27},\\
2^{12} 3^7 J_9,\quad
2^{17} 3^{10} J_{12},\quad
2^{23} 3^{12} J_{15},\quad
2^{27} 3^{15} J_{18},\quad
2^{31} 3^{18} I_{21},\quad
2^{33} 3^{16} J_{21}.
\end{array}
$$
We refer to these normalizations of the Dixmier--Ohno invariants 
as the {\it integral} Dixmier--Ohno invariants (as opposed to the 
{\it classically integral} invariants).  {\it Hereafter we will 
make these normalizations and write $I_{3k}$ or $J_{3l}$ to denote 
the above integral invariants.}

In what follows we invert $I_3$ in order to define six algebraically 
independent functions on the moduli space of quartic plane curves
$$
(i_1,i_2,i_3,i_4,i_5,i_6) = \left( 
\frac{I_6}{I_3^2},
\frac{I_9}{I_3^3},
\frac{I_{12}}{I_3^4},
\frac{I_{15}}{I_3^5},
\frac{I_{18}}{I_3^6},
\frac{I_{27}}{I_3^9}
\right),
$$
and those defining an algebraic extension of $\CC(i_1,\dots,i_6)$:
$$
(j_1,j_2,j_3,j_4,j_5,j_6) = \left( 
\frac{J_9}{I_3^3},
\frac{J_{12}}{I_3^4},
\frac{J_{15}}{I_3^5},
\frac{J_{18}}{I_3^6},
\frac{I_{21}}{I_3^7},
\frac{J_{21}}{I_3^7}
\right)\cdot
$$


\section{Vermeulen Stratification}

In 1983, Vermeulen~\cite{Vermeulen1983} constructed a stratification 
of the moduli space $\Me_3$ of curves of genus 3 in terms of the number 
of hyperflexes and their geometric configuration.  
Similar results were obtained independently around the same time by 
Lugert~\cite{Lugert81}.
The classification of curves by the structure of their Weierstrass points 
identifies more subtle structure of the curves than that provided by 
the automorphism group.

Let $\Me_3^{\circ}$ be  $\{[C] \in \Me_3$, $C$ non-hyperelliptic\!~$\}$, 
$M_s= \{[C] \in \Me_3^{\circ}$, $C$ has at least $s$ hyperflexes\!~$\}$, 
and $M_s^{\circ}= \{[C] \in \Me_3^{\circ}$, $C$ has exactly $s$ 
hyperflexes\!~$\}$.
All strata but $\Me_3^{\circ}$ are closed irreducible subvarieties of
$\Me_3$. Each $M_s^{\circ}$ is the union of the strata with $s$ 
hyperflexes. 
For instance, $M_{12}^{\circ}$ consists of the two moduli points 
corresponding to the Fermat curve and the curve 
$X^4+Y^4+Z^4+3(X^2Z^2+X^2Y^2+Y^2Z^2)=0$ and $M_{11}^{\circ}$ 
and $M_{10}^{\circ}$ are both empty.
In the diagram below, we summarize the relevant data from Vermeulen's
stratification.  We denote by $s$ the number of hyperflexes.

\begin{table}[htbp]
  \centering
\vspace{-4mm}
\begin{minipage}[t]{45mm}
$$
\begin{array}{*{5}c}
\mathcal{X} & & s & \text{dim} & \text{Substrata} \\ 
\hline
\Me_3^\circ 
    & & 0 
    & 6
    & M_1^\circ \\
M_1^\circ 
    & & 1 
    & 5
    & M_2^\circ \\
M_2^\circ 
    & & 2 
    & 4
    & X_1, X_2, X_3 \\
X_2 
    & & 3 
    & 3
    & Y_1, Y_2, Y_3 \\
X_3 & & 3 
    & 3
    & Y_1, Y_3, Y_4, Y_5 \\
X_1 & & 4 
    & 3
    & Y_1 \\
Y_2 & & 4 
    & 2
    & Z_1, Z_5 \\
Y_3 & & 4
    & 2
    & Z_i, 1 \le i \le 8 \\
Y_4 & & 4 
    & 2
    & Z_i,\, i \ne 3, 6, 8 \\ 
Y_5 & & 4
    & 2
    & Z_2, Z_3, Z_6, Z_9 \\
Y_1 & & 5
    & 2
    & Z_1, Z_2, Z_3, Z_4 \\
Z_6 & & 5 
    & 1
    & \Theta, \Pi_i, \Omega_i, \Phi \\
Z_7 & & 5 
    & 1
    & \Pi_i, \Sigma, \Omega_i, \Psi \\
\hline
\end{array}
$$
\end{minipage}
\begin{minipage}[t]{45mm}
$$
\begin{array}{*{5}c}
\mathcal{X} & & s & \text{dim} & \text{Substrata} \\ 
\hline
Z_8 & & 5 
    & 1
    & \Theta, \Pi_i, \Sigma, \Psi \\
Z_2 & &6 
    & 1
    & \Pi_i, \Omega_i, \Phi, \Psi \\
Z_3 & & 6 
    & 1
    & \Theta, \Pi_i, \Omega_i, \Psi \\
Z_5 & &6 
    & 1
    & \Sigma, \Phi, \Psi \\
Z_9 & & 6 
    & 1
    & \Omega_i, \Phi, \Psi \\
Z_4 & & 7 
    & 1
    & \Omega_i, \Psi \\
\Theta  & & 7 
    & 0
    & \\
\Pi_i  & & 7
    & 0
    & \\
Z_1 & & 8 
    & 1
    & \Phi, \Psi \\
\Sigma  & & 8 
    & 0
    & \\
\Omega  & & 9 
    & 0
    & \\
\Phi  & & 12 
    & 0
    & \\
\Psi  & & 12 
    & 0
    & \\

\hline
\end{array}
$$
\end{minipage}
\smallskip 
 \caption{Vermeulen's stratification of $\Me_3$}
  \label{tab:vermeulen}
\vspace{-4mm}
\end{table}

\noindent
{\bf N.B.} The $X_i$ have dimension~$3$, the $Y_i$ have dimension~$2$,
the $Z_i$ have dimension~$1$ and the strata with a Greek letter have
dimension zero.

%


\section{Special Strata of Plane Quartics}

For the special strata $S$ of $\Me_3$ with more than six hyperflexes 
we determine a parametrization of the stratum and a model for a generic 
curve $C/\tilde{S}$ for some finite cover $\tilde{S} \rightarrow S$.
In each case the structure of $\tilde{S} \rightarrow S$ is a Galois 
cover over which the hyperflexes locus splits completely. 
The Dixmier--Ohno invariants are computed over $\tilde{S}$ by their  
sequences of covariants and contravariants, rather than evaluation 
of symbolic expressions.  

\subsection*{Stratum $Z_4$.}

The moduli space $Z_4$ is a one-dimensional subspace of $\Me_3$, for 
which we can find a generic curve defined over $\QQ(i,t)$ of the form:
$$
C : t(t+i)(X^2-YZ)^2 - YZ(2X-Y-Z)(((i-1)t^2+2t+(i+1))X-(it+1)(Y+Z))
$$
where $i^2=-1$ and $t$ is a parameter.  Let $\tilde{Z}_4$ be the 
rational curve with function field $\QQ(i,t)$.  This model parametrizes 
the curve plus the triple of Weierstrass points 
$$
\{(0:0:1),(0:1:0),(1:1:1)\},
$$
with tangent lines defined by $Y=0$, $Z=0$, and $2X-Y-Z=0$, respectively,
and moreover, the remaining hyperflexes of $C$ split over $\QQ(i,t)$.
Thus $\tilde{Z}_4$ defines a pointed moduli space parametrizing three 
hyperflexes, and should determine a Galois cover of the moduli space $Z_4$ 
on which the hyperflex locus splits.

A computation of the Dixmier--Ohno invariants reveils an automorphism 
$t \mapsto i(t-1)/(t+1)$ of order $3$ under which the invariants are 
stable.  The degree 3 quotient given by 
$$
t \mapsto z = \frac{(-(i+1)t^3 + 3t^2 - i)}{(t^2 + (-i+1)t - i)},
$$
provides the Galois cover $\tilde{Z}_4 \rightarrow Z_4$.

The first of the Dixmier--Ohno invariants $(i_1,i_2)$ can be expressed in 
terms of the invariant $z$ as:
$$
\left(
    \frac{(z - 3)(5z - 3)}{4z^2}, 
    \frac{(79z^5-1059z^4-2670z^3+12366z^2-13203z+4455)}{16(z^2-12z+9)z^3}
\right).
$$
Reciprocally, the expression
\begin{equation}
\label{RationalExpressionZ4}
z = \frac{-165(880 i_1^3 + 1336 i_1^2 - 558 i_1 - 1047)}{%
(53680 i_1^3 - 10560 i_1^2 i_2 + 10120 i_1^2 - 57750 i_1 + 7920 i_2 - 6105)}
\end{equation}
gives $z$ as a rational function in $(i_1,i_2)$ so $z$ generates the 
function field of $Z_4$.  We note in particular that $Z_4$ is defined over 
$\QQ$.  

Solving for the algebraic relations in $(i_1, i_2)$ and renormalizing, 
we find a weighted projective equation for $Z_4$ in terms of the first 
Dixmier--Ohno invariants:
\begin{multline}
\label{DefiningPolynomialZ4}
193600 I_6^5 - 35776 I_3^2 I_6^4 + 86784 I_3 I_6^3 I_9 - 961040 I_3^4 I_6^3 \\
  - 2304 I_6^2 I_9^2 + 100608 I_3^3 I_6^2 I_9 - 526608 I_3^6 I_6^2 
  - 65376 I_3^5 I_6 I_9 \\ + 721521 I_3^8 I_6 
  + 1728 I_3^4 I_9^2 - 78048 I_3^7 I_9 + 515889 I_3^{10} = 0.
\end{multline}
Thus from the invariants $I_3$, $I_6$, and $I_9$ we determine a necessary 
condition for a given quartic curve to be in the stratum $Z_4$.  The remaining 
invariants have rational expressions in $z$, so can be readily computed from 
the rational expression~\eqref{RationalExpressionZ4}.  A comparison with the 
remaining Dixmier--Ohno invariants verifies or contradicts the hypothesis that 
a curve with invariants satisfying~\eqref{DefiningPolynomialZ4} lies in $Z_4$.

Given the invariant $z$ for a point in $Z_4$, we can determine a field of 
definition for a representative curve with model $C$ above, by solving for 
a root $t$ of the degree six polynomial
$$
2T^6 + 2(z - 3)T^5 + (z-3)^2T^4 + 2(z^2 - 4z + 1) T^3 
   + 2 z^2 T^2 + 2z(z - 1) T + (z-1)^2,
$$
which is reducible over any field containing a square root of $-1$.


\subsection*{Stratum $Z_1$.}

The moduli space $Z_1$ is a one-dimensional subspace of $\Me_3$ which 
consists of moduli of curves with $8$ hyperflexes and automorphism group 
$D_4$.  There exists a generic curve over some cover $\tilde{Z}_1$ on 
which the of $Z_1$, defined by
$$
C : (t^2+1)(X^2-YZ)^2 = YZ(2X-Y-Z)(2tX-Y-t^2Z).
$$
By computing the Dixmier--Ohno invariants, we find that there exists a 
cyclic degree 4 cover $\tilde{Z}_1 \rightarrow Z_1$.  In particular the 
transformation $t \mapsto (1+it)/(t+i)$ of order four maps $C$ to an 
isomorphic curve.
As above, we find an invariant function $z = t + 1/t + u + 1/u - 1/2$, 
where $u = (1+it)/(t+i)$, such that the absolute Dixmier--Ohno invariants 
can be expressed in $z$.  Specifically, the first two are:
$$
(i_1,i_2) = 
\left(
  \frac{(2z-9)(2z+9)}{4z^2},
  \frac{(2z+9)(8z^2-24z+459)}{2^4z^3}
\right),
$$
Reciprocally, we find an expression for $z$ in terms of the first 
invariants
$$
z = (-153 i_1 - 171)/(26 i_1 - 12 i_2 + 38),
$$
so that $z$ generates the function field of $Z_1$.  The remaining 
invariants can be expressed in terms of the invariant $z$:
\begin{multline*}
(i_3,i_4,i_5,i_6) = 
\left(
  \frac{(2z+9)(4z-27)(8z+9)^2}{2^7z^4},
  \frac{(z+9)^2(2z-45)(2z+9)^2}{4z^5},
  \right.\\
  \left.
  \frac{(z+9)(2z+9)^2(8z+9)(8z^2-129z+837)}{2^5z^6},
  -\frac{(2z - 7)^3}{2^{17}z^9}\right).
\end{multline*}
including $(j_1,\dots,j_6)$:
\begin{multline*}
\\ \Big(
\frac{(2 z + 9) (2 z^2 - 11 z - 9)}{4 z^3},
\frac{(2 z + 9) (4 z^3 - 16 z^2 - 99 z - 1215)}{2^3 z^4},
\frac{(2 z + 9)^2(z - 6)^2}{2 z^4},\\
\frac{(2 z + 9)^2 (z - 6) (8 z + 9) (8 z^2 - 49 z - 18)}{2^5 z^6},
\frac{(2 z + 9)^3 (2 z - 15)^2 (2 z^2 - 3 z + 27)}{2^4 z^7},\\
\frac{(2 z + 9)^2 (56 z^5 - 748 z^4 + 1122 z^3 + 20907 z^2 
- 38880 z - 374706)}{16 z^7}\Big).
\end{multline*}
As in the case of $Z_4$, the moduli space $Z_1$ is defined over $\QQ$.
Given the invarinant $z$ for a point in $Z_1$, we can find a curve $C$ 
which is defined in terms of a root $t$ of the degree four polynomial:
$$
(2 T^4 - T^3 + 12 T^2 - T + 2) - 2 z (T^3 + T),
$$
defining a cyclic cover of degree 4 over $Z_1$.

\section{Zero Dimensional Strata}

It remains to classify the strata of dimension zero in terms of moduli, 
which we denote $\Theta$, $\Pi_1$, $\Pi_2$, $\Sigma$, $\Omega_1$, 
$\Omega_2$, $\Phi$, and $\Psi$.
In the case of $Z_1$ and $Z_4$, the known models for curves in these 
families were not obviously definable over their field of moduli.
With the exception of $\Theta$ and $\Sigma$ below, we will see that 
these exceptional strata have a curve which may be defined over its 
field of moduli.

\subsection*{Stratum $\Theta$.}

The moduli space $\Theta$ is a zero-dimensional subspace of $\Me_3$ 
which is represented by the curve $C_\Theta$ of the form
$$
637(X^2-YZ)^2 = (2X-Y-Z)(aX+bY+cZ)YZ
$$
where 
$$
\begin{array}{l}
a = (-132i-240)t^2+(702i+1068)t-219i+59\\
b = (-88i-160)t^2+(468i+712)t-146i-810\\
c = (108i+428)t^2+(-806i-2032)t+295i+415
\end{array}
$$
with $i^2 = -1$ and $t^3=((10+i)t/2-1)(t-1)$, represented by the 
Dixmier invariants
$$
(i_1,i_2,i_3,i_4,i_5,i_6) = 
\left(
  -\frac{3^3}{2^2 7^2},
  \frac{3^3 173}{2^4 7^2},
  -\frac{3^5 149}{2^7 7^2},
  -\frac{3^5 10223}{2^2 7^4},
  \frac{3^5 5^2 3527}{2^5 7^4},
  \frac{11^4 13^2}{2^{11} 3^{18}7^5}
\right).
$$
and the point in $\Me_3$ is completely determined by the $i_k$ and the 
additional invariants of Ohno:
$$
(j_1,j_2,j_3,j_4,j_5,j_6) = 
\left(
  -\frac{3^3 11}{2^3 7^2},
  -\frac{3^2 4817}{2^4 7^3},
  \frac{3^3}{2\cdot7^2},
  \frac{3^4 535}{2^5 7^3},
  \frac{3^4 55291}{2^7 7^4},
  -\frac{3^2 5486023}{2^5 7^5}
\right).
$$

\subsection*{Strata $\Pi_1, \Pi_2$.}

The strata $\Pi_i$ are zero-dimensional subspaces of $\Me_3$ which 
are represented by the curves $C_{\Pi_i}$ of the form
$$
49(X^2-YZ)^2 = YZ(2X-Y-Z)(aX + b(Y+Z))
$$
where 
$$
\begin{array}{l}
a = (-52i+46)t^2+(49i+25)t-82i-114\\
b = (5i-37)t^2+(-28i+19)t+83i-6, 
\end{array}
$$
with $i^2 = -1$ and $t$ satisfying $2t^3=(-i+1)t^2+4it+(i+1)$.  
The ideal of relations for the absolute Dixmier invariants 
$(i_1,i_2,i_3,i_4,i_5,i_6)$ for $\Pi = \Pi_1 \cup \Pi_2 \cup \Pi_3$, 
is the degree three ideal
\begin{multline*}
(614656\,i_1^3 + 21952\,i_1^2 - 231516\,i_1 - 62613,\\
- 2^5 3^2 2297\,i_2 + 15135904\,i_1^2 - 681236\,i_1 + 418761,\\
- 2^7 2297\,i_3 + 7519344\,i_1^2 - 7084828\,i_1 - 3230271,\\
2^2 2297\,i_4 + 14936160\,i_1^2 + 20448508\,i_1 + 5686083,\\
- 2^7 7^2 2297\,i_5  + 12234260416\,i_1^2 + 10161115868\,i_1 + 1386276669,\\
- 2^{21} 3^{18} 7^3 2297\,i_6+ 87127555902240\,i_1^2 - 19953560617372\,i_1 - 29171717887351)
\end{multline*}
We note, in addition, that the curves $C_{\Pi_i}$ admit an automorphism 
$\sigma(X,Y,Z) = (X,Z,Y)$, which induces a quotient to a non-CM elliptic 
curve $E_{\Pi_i}$, whose $j$-invariant satisfies
$$
j^3 - \frac{3907322953}{3^4 7^4}j^2 
    + \frac{429408710168}{3\,7^4}j 
    - \frac{126488474356752}{7^4} = 0.
$$

\subsection*{Stratum $\Sigma$.}

The stratum $\Sigma$ is represented by the curve 
$$
X^4 + Y^4 + 6\sqrt{-7}X^2Y^2-3(-1+\sqrt{-7})XYZ^2=(7+3\sqrt{-7})/8Z^4.
$$
with absolute Dixmier invariants
$$
(i_1,i_2,i_3,i_4,i_5,i_6) = 
    \Big(
    \frac{3^3}{2^2 7},
    \frac{1557}{2^4 7},
    -\frac{18225}{2^7 7},
    -\frac{28403}{2^2 7},
    \frac{2419065}{2^5 7^2},
    \frac{1}{2^{13} 3^{18} 7^4}\Big),
$$
and absolute Ohno invariants
$$
(j_1,j_2,j_3,j_4,j_5,j_6) = 
\Big(
   -\frac{159}{56},
   -\frac{3249}{112},
   9,
   \frac{14445}{224},
   \frac{166617}{896},
   -\frac{2076561}{1568}
\Big)\cdot
$$


\subsection*{Strata $\Omega_1$, $\Omega_2$.}

These strata consists of two curves 
$$
(X^2-YZ)^2 = (3\pm\sqrt{7})(2X-Y-Z)(X-Y-Z)YZ.
$$
In terms of the absolute Dixmier invariants $(i_1,i_2,i_3,i_4,i_5,i_6)$,
the ideal of relations for $\Omega_1 \cup \Omega_2$ is the degree two 
ideal
\begin{multline*}
(64\,i_1^2 + 64\, i_1 + 9,
    1864\,i_1+64\,i_2-153,
512\,i_3+66416\,i_1+15435,\\
32\,i_4+28504\,i_1+16695,
64\,i_5+383138\,i_1+37737, \\
1624959306694656\,i_6+34973684392\,i_1+5920507885).
\end{multline*}


\subsection*{Stratum $\Phi$.}

The stratum $\Phi$ is represented by the Fermat quartic 
$$
X^4 + Y^4 + Z^4 = 0,
$$
with absolute Dixmier invariants
$$
(i_1,i_2,i_3,i_4,i_5,i_6) = (0,0,0,0,0,-2^{4} 3^{18}),
$$
and absolute Ohno invariants
$$
(j_1,j_2,j_3,j_4,j_5,j_6) = (0,0,0,0,0,0).
$$


\subsection*{Stratum $\Psi$.}

The stratum $\Psi$ is represented by the quartic 
$$
X^4 + Y^4 + Z^4 + 3(X^2 Z^2 + X^2 Y^2 + Y^2 Z^2)
$$
with absolute Dixmier invariants 
$$
(i_1,i_2,i_3,i_4,i_5,i_6) = 
  \Big(\frac{9}{16},  \frac{3^3 7^2}{2^7}
     \frac{3^4 7^3}{2^{11}},\frac{3^5 7^3}{2^7}, \frac{3^6 7^4}{2^11},
     -\frac{5^6}{2^{26} 3^{18}} \Big), 
$$
and absolute Ohno invariants
$$
(j_1,j_2,j_3,j_4,j_5,j_6) = 
  \Big(\frac{63}{2^5}, \frac{2457}{2^7}, \frac{9}{2}, 
       \frac{3969}{2^7}, \frac{177957}{2^{12}}, \frac{606879}{2^{10}}
  \Big)\cdot
$$

\end{document}